\documentclass{article}

\usepackage{latexsym}
\usepackage{amsfonts}
\usepackage{amssymb}
\usepackage{amsmath}

\newtheorem{Theorem}{Theorem}[section]

\newtheorem{Problem}{Problem}[section]
\newtheorem{Conjecture}{Conjecture}[section]

\newcommand{\vc}[1]{\mathbf{#1}}
\newcommand{\C}{\mathcal{C}}
\newcommand{\A}{\mathcal{A}}

\title{Heuristics for The Whitehead Minimization Problem}
\author{R.M. Haralick, A.D. Miasnikov and A.G. Myasnikov   }

\begin{document}

\maketitle

\begin{abstract}
In this paper we discuss several heuristic strategies which allow
one to  solve the Whitehead's minimization problem much faster (on
most inputs) than the classical Whitehead algorithm. The mere fact
that these strategies work in practice leads to several
interesting mathematical  conjectures. In particular, we
conjecture that the length of most non-minimal elements in a free
group can be reduced  by a  Nielsen automorphism which can be
identified by inspecting the structure of the corresponding
Whitehead Graph.
\end{abstract}

\section{Introduction to Whitehead method}

\label{sub-sec:WminProb}
\label{sec:WDA}

Let $X = \{x_1, \ldots, x_m\}$ be a finite alphabet, $X^{-1} =
\{x^{-1} \mid x \in X\}$ be the set of formal inverses of letters
from $X$ and $X^{\pm 1} = X \cup X^{-1}$. A word $w = y_1 \ldots
y_n $ in the alphabet $X^{\pm 1}$ is called {\em reduced} if $y_i
\neq y_{i+1}$ for $i = 1, \ldots, n-1$ (here we assume that
$(x^{-1})^{-1} = x$).    Applying reduction rules $xx^{-1}
\rightarrow \varepsilon, x^{-1}x \rightarrow \varepsilon$ (where
$\varepsilon$ is the empty word), one can reduce each word $w$ in
the alphabet $X^{\pm 1}$ to a reduced word $\overline {w}$. The
word $\overline {w}$ is uniquely defined and does not depend on a
particular sequence of reductions. Denote by $F = F(X)$ the set of
reduced words over $X^{\pm 1}$. The set $F$ forms  a group with
respect to the multiplication $u \cdot v = \overline{uv}$, which
is called   a {\em free } group with {\em basis} $X$. The
cardinality $|X|$  is called the {\em rank}
 of $F(X)$. Sometimes we write $F_n$ instead of $F$ to indicate that
 the rank of $F$ is equal to $n$.

A bijection $\phi: F \rightarrow F$ is called an {\em
automorphism} of $F$ if $\phi(uv) = \phi(u)\phi(v)$ for every $u,
v \in F$. The set  $Aut(F)$ of all automorphisms of $F$ forms a
group  with respect to  composition of maps. Every automorphism
$\phi \in Aut(F)$ is completely determined by the images $\phi(x)$
of elements  $x \in X$. The following two subsets of $Aut(F)$ play
an important part in group theory and topology.

 An automorphism  $t \in Aut(F)$  is called a
 {\em Nielsen automorphism} if  for some $x \in X$ $t$ fixes
 all elements $y \in X, y \neq x$ and maps $x$ to one of the
 elements $x^{-1}$, $y^{\pm 1}x$, $xy^{\pm 1}$. By $N(X)$ we
 denote the set of all Nielsen automorphisms of $F$.

 An automorphism  $t \in Aut(F)$  is called a
 {\em Whitehead automorphism} if either $t$ permutes elements of $X^{\pm 1}$ or
  $t$ fixes a given element $a \in X^{\pm 1}$ and maps each element $x
\in X^{\pm 1}, x \neq a^{\pm 1}$ to one of the elements $x$, $xa$,
$a^{-1} x$, or $a^{-1} x a$. Obviously, every Nielsen
automorphism is also a Whitehead automorphism.  By $W(X)$ we
denote the set of non-trivial  Whitehead's automorphisms of the
second type.

Observe that
  $$|N(X)| = 4n(n-1), \ \ \ |W(X)| = 2n4^{(n-1)} - 2n$$
where $n = |X|$ is the rank of $F$.

It is known \cite{Lyndon77compbgt} that every automorphism from $Aut(F)$ is a
product of finitely many Nielsen (hence Whitehead)  automorphisms.

 The automorphic orbit
$Orb(w)$ of  a  word $w \in F$ is the set of all automorphic
images of $w$ in $F$:
\[Orb(w) = \{ v \in F  \mid  \exists \varphi \in Aut(F) \mathrm{\; such \;that}\; w^\varphi =
v\}.\] A word $w \in F$ is called {\em minimal} (or {\em
automorphically minimal}) if $|w| \leq |w^\varphi|$ for any
$\varphi \in Aut(F)$. By $w_{min}$ we denote a word of minimal
length in $Orb(w)$. Notice that $w_{min}$ is not  unique.

\begin{Problem}[Minimization Problem (MP)]
\label{RP} For a word $u \in F$ find an automorphism $\varphi \in
Aut(F)$ such that $u\varphi = u_{min}$.
\end{Problem}

 In 1936  J. H. C. Whitehead  proved the following result which gives
 a solution to the minimization problem  \cite{Whitehead36equivsets}.

\begin{Theorem}[Whitehead]
\label{thm:T1} Let $u,v \in F_n(X)$ and $v \in Orb(u)$. If
$|u|>|v|$, then there exists $t \in W(X)$ such that
\[
|u|>|ut|.
\]
\end{Theorem}

An automorphism $\phi \in Aut(F)$  is called {a length-reducing}
automorphism for a given  word $u \in F$ if $|u\phi| < |u|$. The
theorem  above claims that the finite set $W(X)$ contains a
length-reducing automorphism for every non-minimal word $u \in F$.
This allows one to design a simple search algorithm for (MP).

 Let $u \in F$. For each $t \in W(X)$ compute the length
of the tuple $ut$ until $|u|
> |ut|$, then put $t_1 = t, u_1 = ut_1$. Otherwise stop and output
$u_{min} = u$. The procedure above is called the {\em Whitehead
Length Reduction}  routine (WLR). Now  Whitehead Reduction
Algorithm (WRA) proceeds as follows. Repeat WLR on $u$, and then
on $u_1$, and so on, until on some  step $k$ WRL gives an output
$u_{min}$. Then $ut_1 \ldots t_{k-1} = u_{min}$, so $\phi = t_1
\ldots t_{k-1}$ is a required automorphism.

 Notice,  that the iteration procedure WRA  simulates the classical
 greedy  descent method ($t_1$ is a successfull  direction from $u$,
 $t_2$ is a successfull direction from $u_1$, and etc.). Theorem \ref{thm:T1}
guarantees that the greedy approach will always converge to the
global minimum.

Clearly, there could be at most $|u|$ repetitions of WLR on an
input $u \in F$
\[|u| > |ut_1|>...>|ut_1 ... t_l| = u_{min}, \ \ \ l \leq |u|.\]
 Hence the worst case complexity of the algorithm WRA is bounded from above by
\[
cA_n|u|^2,
\]
where $A_n = 2n4^{(n-1)} - 2n$ is the number of Whitehead
automorphisms in $W(X)$. Therefore, in the worst case scenario,
the algorithm seems to be impractical for free groups with large
ranks. One can try to improve on the number of steps which takes
to find a length-reducing automorphism for a given non-minimal
element from $F$. In this context the main question of interest is
the complexity of the following

\begin{Problem}[Length Reduction Problem]
For a given non-minimal element $u \in F$ find a length-reducing automorphism.
\end{Problem} We refer to \cite{ MM03wga} for a general discussion on this problem.

 In the next section we  give some empirical
evidence that using  smart strategies in selecting Whitehead
automorphisms $t \in W(X)$  one can dramatically improve the
average complexity of WRA in terms of the rank of a group.

\section{Heuristics for Length Reduction Problem}
\label{sec:Strategies}

\subsection{Nielsen first}

 The first heuristic comes from a very naive approach: replace $W(X)$
 by $N(X)$ in the Whitehead length reduction routine WLR  and denote the resulting
  routine by NLR. Since the size of $N(X)$ is quadratic and the size of $W(X)$ is
exponential in  the rank of $F$, the algorithm  NRA may give a
real speedup in computations.
 However,  it is known (see \cite{Lyndon77compbgt}) that the Whitehead theorem
 above does not hold after replacement of
$W(X)$ by $N(X)$. Therefore, the algorithm NRA will not give the
correct answer at least on some inputs. But this is not the end of
the story. Now the question is  {\em how often the length
reduction routine NLR gives the correct answer?}

To get some insights, we perform a simple experiment. For free
groups $F_3$, $F_4$ and $F_5$ we generate test sets of non-minimal
elements of Whitehead Complexity 1 (see definitions in
\cite{MM03wga}), described in Table \ref{tab:data_non-min}. For a
detailed description of the data generation procedures  we refer to \cite{HMM03pr}.

\begin{table}[ht]
\begin{center}
\begin{tabular}{|c|c|c|c|c|c|}
\hline
Dataset &  Group  & Dataset Size  &  Min. length & Avg. length &  Max. length \\
\hline
$D_3$ & $F_3$ & 10143 & 3 & 558.2& 1306\\
\hline
$D_4$ & $F_4$ & 10176 & 4 & 570.9 & 1366\\
\hline
$D_5$ & $F_5$ & 10165 & 5 & 581.3& 1388\\
\hline
\end{tabular}
\end{center}
\caption{Statistics of sets of non-minimal elements. }
\label{tab:data_non-min}
\end{table}

For each set $D_n$ we compute the fraction of elements from $D_n$
which have length-reducing Nielsen automorphisms. The results of
the computations together with the corresponding 95\% confidence
intervals are given in Table \ref{tab:frac_non-min}. We can see
that most of the words have been reduced by Nielsen automorphisms.
We would like to mention here, that it can be shown statistically
that increasing the length of elements in the datasets does not
significantly change the results of experiments.

%
%
\begin{table}[ht]
\begin{center}
\begin{tabular}{|l|c|c|c|}
\hline
Dataset & $D_3$ & $D_4$ & $D_5$\\
\hline
Fraction & 0.998 & 0.997  & 0.998 \\
\hline
95\% Conf. Interval & [0.9970,0.9988] & [0.9957,0.9979] & [0.9970,0.9988] \\
\hline
\end{tabular}
\end{center}
\caption{Fraction of  elements in the sets $D_3$, $D_4$ and $D_5$
with length-reducing Nielsen automorphisms.}
\label{tab:frac_non-min}
\end{table}
%
%
Based on these experiments one can speculate that with very
high probability Nielsen automorphisms  reduce the length of a given
non-minimal element in $F$. More precisely, we state the following
\begin{Conjecture}
\label{conj:NielsenFirst}
 Let $U_n$ be the set of all non-minimal elements in $F$ of length $n$
and $NU_n \subset U_n$ the subset of elements  which have Nielsen
length-reducing automorphisms. Then
 $$
\lim_{n \rightarrow \infty} \frac{|NU_n|}{|U_n|}  = 1.
$$
\end{Conjecture}

Our first heuristic is based on this conjecture and simply
suggests to try Nielsen automorphisms first in the routine WLR,
i.e., in this case we assume that in the fixed listing of
automorphisms of $W(X)$ the automorphism from $N(X)$ come first.
We refer to this heuristic as to \emph{Nielsen First} and denote
the corresponding Length Reduction Routine and the Whitehead
Reduction algorithm (with respect to this ordering of $W(X)$) by
$\mathrm{WLR}_{NF}$ and $\mathrm{WRA}_{NF}$.

The expected value of the number of steps for the routine
$\mathrm{WLR}_{NF}$ to find a length-reducing automorphism on an
input $u \in F$ of length $n$ is equal to
\[
P_n|N(X)| + (1-P_n)(|W(X)| - |N(X)|),
\]
where $P_n =  NU_n / U_n$.

Given that Conjecture \ref{conj:NielsenFirst} is true, we expect
$\mathrm{WRA}_{NF}$ to perform much better on average. In the next
section we describe experimental results supporting this strategy.

\subsection{Cluster analysis}

According to the heuristic NF one has to apply  Nielsen automorphisms to a given input
$w$ in some fixed order,
which is independent of the word $w$.  Intuitively, we expect some
automorphisms to be more likely to reduce the length of a given
word than the others. It suggests that
the conditional probabilities
\[
Prob( |w t| < |w| \mid w), \ t \in N(X)
\]
may  not be equal for  different non-minimal words $w \in F$, so
the order in which  Nielsen automorphisms are applied to an input $w$ should depend on the
word $w$ itself.

The question we would like to address next is whether it is
possible to find a dependence between a non-minimal word $w$ and its
length-reducing Nielsen automorphisms. For this purpose we employ
methods from \emph{Statistical Pattern Recognition}.

Briefly, Pattern Recognition aims to classify  a variety of given objects into
 categories based on the existing  statistical information.
The  objects are typically presented by collections of
 measurements or observations (called \emph{features}) which are real numbers.
  In this event the tuple of features that corresponds
to a given object  is called  a \emph{feature vector}, it can be
viewed as  a point in the appropriate multidimensional vector
space $\mathbb{R}^d$. Most of the approaches in Statistical
pattern recognition are based on statistical characterizations of
features, assuming that objects are generated by a probabilistic
system. The detailed description of Pattern Recognition methods is
out of scope of this report. We refer interested readers  to
\cite{Duda00patternclass,Fukunaga90PR,theodoridis99pattern}
 for general  introduction to the subject, and \cite{HMM03pr} for
 applications of pattern recognition methods in groups.

Unsupervised learning or \emph{clustering}  methods of pattern recognition are used when no a priori
information about the objects is available.
In this case there are general algorithms to group  the feature vectors of objects into some "natural classes"
(called \emph{clusters})  relative to the specified similarity assumptions. Intuitively, the objects whose
feature vectors belong to the same cluster
are more similar to each other than the objects with the feature vectors in  different clusters.

The most simple and widely used clustering scheme is called {\em
$K$-means}. It is an iterative method. Let $D = \{\vc{x}_1,
\ldots, \vc{x}_N\}$ be a set of given objects, represented by the
corresponding feature vectors $\vc{x} \in \mathbb{R}^d$. $K$-means
begins with a set of $K$ randomly chosen cluster centers
$\mu^0_1,\ldots, \mu^0_K \in \mathbb{R}^d$. At iteration $i$ each
feature vector is assigned to the nearest cluster center  (in some
metric $||\  ||$ on $\mathbb{R}^d$). This forms the cluster sets
$C^i_1,\ldots, C^i_K$, where

$$
C^i_j = \{ \vc{x}  \mid  ||\vc{x} - \mu^i_j|| \le ||\vc{x} - \mu^i_m||,\ \vc{x} \in D, \ m=1,\ldots,K \}.
$$

Then each cluster center is redefined as the mean of the feature
vectors assigned to the cluster:

$$
\mu^{i+1}_k = \frac{1}{\ |C^i_k| } \sum_{\vc{x} \in C^i_k} \vc{x}.
$$

Each iteration reduces the criterion function $J^i$ defined as
$$
J^i = \sum_{k=1}^K \sum_{\vc{x}\in C^i_k} ||\vc{x} - \mu^i_k||.
$$

As this criterion function is bounded below by zero, the
iterations must converge. This method works well when clusters are
mutually exclusive and compact around their center means.

Here we claim that  $K$-means algorithm allows one to discover some natural classes of non-minimal words.
We show below that analysis of the corresponding cluster structures sheds some light on the relation
between non-minimal words and their length-reducing automorphisms.

We define  features of
elements $w \in F(X)$ as follows.  Recall that the Labelled Whitehead Graph
$WG(w) = (V,E)$ of an element $w \in F(X)$
 is a weighted non-oriented graph, where the set of vertices $V$ is equal to the set  $X^{\pm
 1}$, and for $x_i, x_j \in X^{\pm 1}$ there is an edge $(x_i,x_j) \in E$
if the subword $x_i x_j^{-1}$ (or $x_jx_i^{-1}$) occurs in the word $w$ viewed as a cyclic word.
Every edge $(x_i,x_j)$ is
assigned a weight $l_{ij}$ which is the number of times the
subwords $x_ix_j^{-1}$ and $x_jx_i^{-1}$ occur in $w$.

Let $l(w)$ be a vector of edge weights in the Whitehead Graph $WG(w)$
with respect to a fixed order.
We define a feature vector $f(w)$ by
\[
f(w) = \frac{1}{|w|}l(w).
\]

To execute the $K$-means algorithm one has to define in advance
the expected number of clusters $K$. Since we would like  these
clusters to be related to the set of Nielsen automorphisms $N(X)$
we put $K = | N(X)|$.

To evaluate usefulness of the clustering we use the goodness
measure $R_{max}$ defined below.  Let  $\C \subset D$ be a cluster
of the data set $D \subset F_n = F(X)$. For $t \in N(X)$ define
$$
R(t,\C) = \frac{| \{ w \in \C  \mid |w t |<| w | \}}{|\C|}.$$
 The number
$R(t,\C)$ shows how many elements in $\C$ are
reducible by $t$. Now put
$$
R_{max}(\C) = \max \{R(t,\C) \mid  t \in N(X) \}
$$
and denote by $t_{\C}$ a Nielsen automorphism $t \in N(X)$ such
that $ R(t_\C,\C) = R_{max}(\C)$. The number $R_{max}(\C)$ shows
how many elements in $\C$ can be reduced by a single automorphism,
in this case by $t_{\C}$. We also define the average value of the
goodness measure
\[
avg(R_{max}) = \frac{1}{K}\sum_{i=1}^K R_{max}(\C_i),
\]
where $K$ is the number of clusters.

 The results of $K$-mean cluster analysis of sets of
randomly generated non-minimal elements in free groups $F_3$,
$F_4$, $F_5$ are given in Table \ref{tab:clust_eval}. It shows
that  more that 70\% of elements in every cluster can be reduced
by the same Nielsen automorphism. In the free group $F_3$, where
the number of clusters is significantly smaller, the corresponding
percentage is over 98\%. Moreover, our experiments show  that
$t_{\C_i} \neq t_{C_j}$ for $i\neq j$. In other words there are no
two distinct clusters such that one and the same Nielsen
automorphism reduces most of the elements in both clusters.

%
%
\begin{table}[ht]
\begin{center}
\begin{tabular}{|l|c|c|c|c|}
\hline
Free group  & $F_3$ & $F_4$ & $F_5$ \\
\hline
number of clusters, $K$& 24 & 48 & 80 \\
\hline
$avg(R_{max})$, $K$-means & 0.985  & 0.879& 0.731 \\
 \hline
\end{tabular}
\caption{ Average values of the goodness measure $R_{max}$ for
$K$-means clustering. } \label{tab:clust_eval}
\end{center}
\end{table}

The discovered cluster structure gives rise to the following strategy
in solving the Length Reduction Problem for a given word $w$.
Let $\mu_1, \ldots, \mu_K$ be the centers of clusters $\C_1, \ldots, \C_K$ computed by the
$K$-means procedure. We compute the distance $||f(w) - \mu_i||$ for each $i = 1, \ldots, K$.
Now we list the Nielsen automorphisms in $N(X)$ in the order
 $t_{i_1}, t_{i_2}, \ldots, t_{i_K}$ with respect to the distances
\[
||f(w) - \mu_{i_1}|| \leq ||f(w) - \mu_{i_2}|| \leq \ldots \leq  ||f(w) - \mu_{i_K}||.
\]

To find a length reducing automorphism for a given word $w$ we subsequently apply
automorphisms from $N(X)$ in the prescribed order until we find an automorphism
$t_i \in N(X)$ which reduces the length  of $w$. If such an automorphism does not exist
 we proceed with the remaining automorphisms from $W(X) - N(X)$ as in the  NF heuristic.

 From the description of the $K$-means method we know that
clusters are characterized by the center means of the feature
vectors of elements in the same cluster. The observations above
lead us to the following vaguely stated conjecture, which gives a
model to describe behavior of non-minimal elements from $F$ in
terms of their feature vectors.

\begin{Conjecture}
\label{conj:hypersurf}
The feature vectors of weights of the Whitehead Graphs of elements
from $F$ are separated into bounded  regions in the corresponding
 space. Each such  region can be bounded by a hypersurface and
corresponds to a particular Nielsen automorphism  in a sense that
all elements in the
 corresponding class can be reduced by that automorphism.
 \end{Conjecture}

\subsection{Improvement on the clustering}

Experiments with $K$-means clustering algorithm show that
clustering is a useful tool in solving the length reduction problem.
Now, the goal is to make clustering more effective.  The further
analysis of the clusters  suggests that to some extent they
correspond to partitions of elements in $F$ which can be reduced
by one and only one Nielsen automorphism. To verify this
conjecture we perform the  following experiment.

  Let $S \subset F_n = F(X)$ be a set of randomly generated non-minimal
  elements   and $D$ the set used for cluster analysis in the previous
  section. Note that $S$ is generated  independently from the set
  $D$.
 For each automorphism $t \in N(X)$ put
 $$
O_t = \{ w \in S \mid  \forall r \in N(X) (|wr| < |w|
\Longleftrightarrow r = t) \}
 $$
 and define new cluster centers by
\begin{equation}
\label{eq:centers}
\lambda_t = \frac{1}{|O_t|} \sum_{w \in O_t} f(w)
\end{equation}
as the mean feature vector of the elements from $S$ that can be
reduced only by $t$ and no other automorphisms.

We cluster elements from $D$ based on the distance between the
corresponding feature vector and centers $\lambda_t$:
\[
\C_t = \{ w \in D \mid  \forall r \in N(X) (||f(w) - \lambda_t|| \le
||f(w) - \lambda_r||)  \}.
\]

The results of evaluation of  the clusters $\C_t$ are given in
Table \ref{tab:clust_eval-2}. One can see that the goodness
measure is improved and is close to 1 in every case.

%
%
\begin{table}[ht]
\begin{center}
\begin{tabular}{|l|c|c|c|c|}
\hline
Free group  & $F_3$ & $F_4$ & $F_5$ \\
\hline
number of clusters, $K$& 24 & 48 & 80 \\
\hline
 $avg(R_{max})$, distance to $\lambda_t$ & 0.998& 0.993& 0.991 \\
 \hline
\end{tabular}
\caption{ Average values of the goodness measure $R_{max}$ for
the clustering based on the distance to
the estimated centers $\lambda_t$. }
\label{tab:clust_eval-2}
\end{center}
\end{table}

Similar to the strategy based on  the centers of the $K$-means
clusters, we define a new search procedure $\mathrm{WRA}_{C}$
which employs a heuristic based on the distances to centers
$\lambda_t$. Let $w$ be a word and $<\lambda_{t_1}, \ldots, \lambda_{t_K}>$ be
the centers  corresponding to each of the Nielsen automorphisms
$t_i \in N(X)$. Put $d(i) = ||f(w) - \lambda_{t_i}||$ and construct a
vector
 $$<d(m_1), d(m_2), \ldots, d(m_K)>,$$
  where
   $$d(m_1) \leq d(m_2) \leq \ldots \leq d(m_K).$$

To find a length reducing automorphism for a given word $w$, the
algorithm $\mathrm{WRA}_{C}$ applies Whitehead automorphisms  to
$w$ in the following order. First, Nielsen automorphisms $t_{m_1},
\ldots, t_{m_K}$ are applied subsequently. If none of the Nielsen
automorphisms reduces the length of $w$ the algorithm
$\mathrm{WRA}_{C}$ proceed with the remaining automorphisms
$W(X) - N(X)$ in some fixed order.

Based on the results of the cluster analysis from Table
\ref{tab:clust_eval-2}, we expect the algorithm $\mathrm{WRA}_C$ to
reduce a non-minimal word $w$ using very few elementary
automorphisms on average.

\subsection{Maximal weight edges}

Now we would like to take a closer look at the edges' weight
distributions at the cluster centers. First, observe that every
edge in the Whitehead graph $WG$, except for the ones which
correspond to subwords of type $x^2$, $x \in X^{\pm 1}$, will
correspond to subwords reducible by two particular Nielsen
transformations. For example, edge connecting nodes $a$ and $b$
corresponds to subwords $(ab^{-1})^{\pm 1}$ both of which are
reduced by automorphisms
\begin{eqnarray*}
(a  \rightarrow  a b,  b    \rightarrow  b), \\
(a  \rightarrow  a,    b   \rightarrow   b a).
\end{eqnarray*}
In fact there is no other Nielsen transformation that will reduce the length of words $(a b^{-1})^{\pm 1}$.

To generalize, let $WG(w)$ be a Whitehead graph of a word $w$ with the vertex set $V$ and the set of edges $E$.
Let $e = (x, y^{-1})$, $x, y^{-1} \in V$, be an edge in $E$. By construction $e$ corresponds
to subwords $s_e = (x y)^{\pm 1}$ of the word $w$. The only Nielsen automorphisms which reduce
length of the subwords $s_e$ are
$$
\psi^x_e: x \rightarrow x y^{-1},\ z \rightarrow z, \   \forall  z \neq x,\ z \in X
$$
and
$$
\psi^y_e: y \rightarrow x^{-1} y,\ z \rightarrow z, \ \forall  z \neq y,\ z \in X.
$$
We will call automorphisms $\psi^x_e, \psi^y_e$ the  length
reducing Nielsen automorphisms with respect to the edge
$e=(x,y^{-1})$ and denote $\psi_e=\{\psi^x_e, \psi^y_e\}$.

The following  phenomenon has been observed for all clusters  in free groups $F_3$, $F_4$, and $F_5$.
Let $\C_t$ be a cluster of a test set $D_n$, $n=3,4,5$, then for all $t \in N(X)$,
\[
t \in \psi_{e_{max}},
\]
where $e_{max}$ is the edge having the maximal weight in the cluster center $\lambda_t$.
It suggests that at least in the case of free groups $F_3$, $F_4$, $F_5$ one can try to estimate
a length-reducing automorphism for given word $w$ by taking the length-reducing
Nielsen automorphisms of the highest weight edge in the Whitehead graph $WG(w)$.

To evaluate the goodness of the heuristic based on the maximal edge weight in
the Whitehead graph we compute the
fraction of elements in the sets $D_3$, $D_4$ and $D_5$, reducible by the Nielsen automorphisms
corresponding to the maximal weight edge. The corresponding goodness
measure, evaluated on a set $D$, is given by
\[
\mathcal{G}_{MAX} = \frac{1}{|D|}|\{ w \in D \mid \exists t \in
\psi_{e_{max}(w)}, \ \mathrm{s.t.}\ |wt| < |w|\}|.
\]

%
\begin{table}[ht]
\begin{center}
\begin{tabular}{|l|c|c|c|c|}
\hline
Dataset  & $D_3$ & $D_4$ & $D_5$ \\
\hline
 $\mathcal{G}_{MAX}$& 0.991  & 0.986 & 0.986  \\
\hline
\end{tabular}
\caption{ Values of the goodness measure $\mathcal{G}_{MAX}$ for sets of non-minimal
elements in free groups $F_3$, $F_4$ and $F_5$. }
\label{tab:max_weight_eval}
\end{center}
\end{table}

Values of the  goodness measure $\mathcal{G}_{MAX}$ for test sets in
free groups $F_3$, $F_4$ and $F_5$
are given in Table \ref{tab:max_weight_eval}.
It shows, that heuristic is surprisingly effective. Nevertheless, centroid based
method still yields better results. Note that  $\mathcal{G}_{MAX}$
measures success of applying two automorphisms corresponding to the maximal
weight edge, where the centroid based  method was evaluated by the success rate
of only one automorphism which corresponds to the closest center.

The observation provides a new search procedure which we denote by
$\mathrm{WRA}_{MAX}$. Let $w$ be a word and $WG(w) = (V,E)$ be the
corresponding Whitehead graph. Denote by $E'$ the set of edges
which do not correspond to the subwords of type $x^{\pm 2}$, $x
\in X$
\[
E' = \{ e \in E \mid e \neq (v,v^{-1}), v \in V\}.
\]

It has been shown above, that for each edge $e$ from $E'$ there exists
 two unique length reducing automorphisms. Note that $2  |E'| = |N(X)|$,
where $N(X)$ is the set of Nielsen automorphisms for free group $F(X)$.

We can order Nielsen automorphisms $\psi_{e_i} \subset N(X)$:
\begin{equation}
\label{eq:max_weight_order}
<\psi_{e_1}, \psi_{e_2}, \ldots, \psi_{e_{|E'|}}>
\end{equation}
such that edges $e_1, \ldots , e_{|E'|}$ are chosen
according to  the decreasing order of the values of the corresponding weights
\[
\omega_{e_1} \geq \omega_{e_2} \geq \ldots \geq \omega_{e_{|E'|}}.
\]
Note that $\psi_e$ is not a single automorphism, but a pair of Nielsen
length reducing automorphisms with respect to the edge $e$. Here we
do not give any preference in ordering automorphisms in $\psi_e$.

To find a length-reducing automorphism for $w$ procedure
$\mathrm{WRA}_{MAX}$ first  applies Nielsen automorphisms in the
order given by (\ref{eq:max_weight_order}). If none of the
Nielsen automorphisms reduces the length of $w$,
$\mathrm{WRA}_{MAX}$ proceeds with the remaining automorphisms
from  $W(X) - N(X)$.

\section{Comparison of the strategies}

In this section we describe experiments designed  to compare the
performance of WRA implemented with different search strategies.
We compare four variations of the algorithm. $\mathrm{WRA}_R$ is
the variation of WRA, where a random order of the elements from
$W_n$ is used when searching for a length reducing
automorphisms. $\mathrm{WRA}_{NF}$ and $\mathrm{WRA}_C$ correspond
to the implementations with  Nielsen First and Centroid based
heuristics respectively. The algorithm $\mathrm{WRA}_{MAX}$
employs strategy which applies automorphisms corresponding to the
largest edge weights of the Whitehead Graph. The algorithms were
compared on randomly generated sets of primitive elements $S_3$, $S_4$, $S_5$  in free
groups $F_3$, $F_4$, and $F_5$, respectively. Some descriptive statistics of the test sets $S_n$
are given in Table \ref{tab:prim_datasets}.

\begin{table}[ht]
\begin{center}
\begin{tabular}{|c|c|c|c|c|c|}
\hline
Dataset &  Group  & Dataset Size  &  Min. length & Avg. length &  Max. length \\
\hline
$S_3$ & $F_3$ & 5645 & 3 & 1422.1 & 143020 \\
\hline
$S_4$ & $F_4$ & 5241 & 4 & 2513.1   & 168353 \\
\hline
$S_5$ & $F_5$ & 3821 & 5 & 2430.5 & 160794 \\
\hline
\end{tabular}
\end{center}
\caption{Statistics of the test sets of primitive elements. }
\label{tab:prim_datasets}
\end{table}

Let ${\mathcal A}$ be one of the  variations $\mathrm{WRA}_R$,
$\mathrm{WRA}_{NF}$, $\mathrm{WRA}_C$, and $\mathrm{WRA}_{MAX}$ of
the Whitehead Reduction Algorithm. By an elementary step of the
algorithm $\A$  we mean one application of a Whitehead
automorphism to a given word. Below we evaluate  the performance
of $\A$  with respect to the number of elementary steps that are
required by $\A$  to execute a particular routine.

Let $N_{total} = N_{total}(\A,S_n)$ be the average of the total
number of elementary steps required by $\A$ to reduce a given
primitive element $w \in S_n$ to a generator.

By $N_{red} = N_{red}(\A,S_n)$ we denote the average number of
elementary length-reducing steps required by $\A$ to reduce a
given primitive element $w \in S_n$ to a generator, so $N_{red}$
is the average number of "productive" steps performed by  $\A$.
 It follows that if $t_1, \ldots, t_l$
are all the length reducing automorphisms found by $\A$ when
executing its routine on an input $w \in S_n$ then $|w t_1 \ldots
t_l| = 1$ and the average value of $l$ is equal to $N_{red}$.

Finally, denote by  $N_{LRP} = N_{LRP}(\A,S_n)$ the average number
of elementary steps  required by $\A$  to find a length-reducing
automorphism for a given non-minimal input $w$.

\begin{table}[ht]
\begin{center}
\begin{tabular}{|l|r|r|r|}
\hline
Strategy      &  $N_{total}$ & $N_{red}$ &   $N_{LRP}$ \\
\hline
$\mathrm{WRA}_C$  & 19.9  & 18.4 & 1.1  \\
\hline
$\mathrm{WRA}_{MAX}$ & 47.1  & 23.9 & 1.9 \\
\hline
$\mathrm{WRA}_{NF}$ & 207.8 & 28.2 & 7.37  \\
\hline
$\mathrm{WRA}_{R}$ & 374.8 & 29.8 & 12.6  \\
\hline
\end{tabular}
\end{center}
\centerline{a) $F_3$;}
\begin{center}
\begin{tabular}{|l|r|r|r|}
\hline
Strategy      &  $N_{total}$ & $N_{red}$ &   $N_{LRP}$\\
\hline
$\mathrm{WRA}_C$ & 58.8  & 34.1 & 1.4  \\
\hline
$\mathrm{WRA}_{MAX}$ &  152.8 & 42.5 & 3.0 \\
\hline
$\mathrm{WRA}_{NF}$ & 1052.6 & 56.2 & 18.7 \\
\hline
$\mathrm{WRA}_{R}$ & 2610.4 & 58.8 & 44.4  \\
\hline
\end{tabular}
\end{center}
\centerline{b) $F_4$;}
\begin{center}
\begin{tabular}{|l|r|r|r|}
\hline
Strategy      &  $N_{total}$ & $N_{red}$ &   $N_{LRP}$ \\
\hline
$\mathrm{WRA}_C$ & 162.0  & 50.9 & 2.4 \\
\hline
$\mathrm{WRA}_{MAX}$ & 342.2  & 58.8 & 4.5 \\
\hline
$\mathrm{WRA}_{NF}$ & 2307.6 & 75.4 & 30.6  \\
\hline
$\mathrm{WRA}_{R}$ & 15939.6  & 78.8 & 202.0  \\
\hline
\end{tabular}
\end{center}
\centerline{c) $F_5$.}

\caption{Results of experiments with sets of primitive elements in
free groups $F_3$,$F_4$ and $F_5$. Counts are averaged over  all
inputs. } \label{tab:res}
\end{table}

 In Table \ref{tab:res} we present results of our experiments on performance
of the algorithms  $\mathrm{WRA}_R$, $\mathrm{WRA}_{NF}$,
$\mathrm{WRA}_C$ and $\mathrm{WRA}_{MAX}$ on the test sets $S_n$,
$n = 3,4,5$.  The algorithms compare as expected. The algorithms
$\mathrm{WRA}_C$ and $\mathrm{WRA}_{MAX}$
 perform very efficiently with the
numbers $N_{LRP}$ and $N_{red}$ being small.  Algorithm
$\mathrm{WRA}_C$  based on  the centroid approach shows best over
all performance. The growth of the numbers $N_{LRP}$  with the
rank could be  explained by occasional occurrence of non-minimal
words that cannot be reduced by Nielsen automorphisms. In this
event the algorithm tries Whitehead automorphisms from $W(X)-
N(X)$ the number of which growth exponentially with the rank.
Notice that in every  our experiment the number $N_{red}$ of
length reductions performed by $\mathrm{WRA}_C$ is less than the
corresponding number in the other approaches.

In Table \ref{tab:corr}  we give the correlation coefficients
showing dependence of
the number of elementary steps required by a particular algorithm
to find a length-reducing automorphism with respect to the length
of the input words.
The coefficients are negative in all cases which indicates that the
values of $N_{LRP}$ do not increase when the words' length increases.

\begin{table}[ht]
\begin{center}
\begin{tabular}{|l|r|r|r|r|}
\hline
Strategy      &  $F_3$ & $F_4$ &   $F_5$ \\
\hline
$\mathrm{WRA}_C$  &  -0.008 & -0.001 &  -0.006\\
\hline
$\mathrm{WRA}_{MAX}$ & -0.024  & -0.023  &  -0.038\\
\hline
$\mathrm{WRA}_{NF}$ & -0.009 &  -0.022 &   -0.022\\
\hline
$\mathrm{WRA}_{R}$ & -0.038 & -0.014 &   -0.035\\
\hline
\end{tabular}
\end{center}
\caption{Correlation coefficients between words length and values of
  $N_{LRP}$. Negative coefficients indicate that $N_{LRP}$ does not
increase when length increases. }
\label{tab:corr}
\end{table}

\section{Conclusions}

The experimental results presented in this paper show that using
appropriate heuristics in the  algorithm WRA, one can significantly reduce
the complexity of the Whitehead minimization problem on most inputs
with respect to the group rank. Suggested heuristics reduce the average number
of Whitehead automorphisms required to find a length-reducing
automorphism for a given word. The performance of heuristic algorithms
tested on the sets of randomly generated primitive elements shows robust
behavior and does not deteriorate when the
length of the input words increases.

One of the interesting  contributions of this paper is  the
empirically discovered properties of non-minimal elements of free
groups formulated in Conjectures \ref{conj:NielsenFirst} and
\ref{conj:hypersurf}. These conjectures suggest that the length of
a "generic" non-minimal elements in a free group can be reduced by
a  Nielsen automorphism. Moreover, the feature vectors of the
weights of the Whitehead's Graphs of non-minimal elements  are
divided into "compact" regions in the corresponding vector  space.
Each such region is related  to a particular Nielsen automorphism,
that reduces the length of all elements in the region.  We believe
this is one of those  few cases  when a meaningful rigorous, but
not intuitively clear conjecture, in group theory was obtained by
using experimental simulations and statistical analysis of the
problem.

It remains to be seen why  the algorithm $\mathrm{WRA}_C$ is able
to find minimal elements using a smaller number of length reductions on average. We are
going to address this issue in the subsequent paper.

%
%
%

\end{document}